\documentstyle[12pt]{amsart}

\def\bC{{\Bbb C}}

\def\bZ{{\Bbb Z}}

\newtheorem{thm}{Theorem}[section]
\newtheorem{lemma}[thm]{Lemma}
\newtheorem{prop}[thm]{Proposition}
\newtheorem{cor}[thm]{Corollary}
\newtheorem{defi}[thm]{Definition}
\newtheorem{rk}[thm]{Remark}

\def\prf{\noindent{\textsc{Proof}}\rm\ }
\def\endprf{\ \hfill $\Box$}
\begin{document}
\baselineskip=16pt
\parindent=0pt

\title[]{Insertion and Elimination Lie Algebra: the ladder case.}

\author[]{Igor Mencattini$^\ddagger$ and Dirk Kreimer$^\dagger$}
{ {\renewcommand{\thefootnote}{\fnsymbol{footnote}}
\footnotetext{\kern-15.3pt AMS Subject Classification: 17B65, 17B70,
16W30, 81T18, 81T15.}
}}

{ {\renewcommand{\thefootnote}{\fnsymbol{footnote}}
\footnotetext{\kern-15.3pt $\dagger$, $\ddagger$ supported in parts by NSF grant 0205977; Ctr. Math. Phys. at Boston U., BUCMP/03-04.}
}}

%\author[]{Igor Mencattini}

\address{Boston University, Department of Mathematics and Statistics,
Boston University, 111 Cummington Street,
Boston, MA 02215, USA}

\email{igorre@@math.bu.edu}

%\author[]{Dirk Kreimer}

\address{CNRS at IHES, 35,
route de Chartres, 91440, Bures-sur-Yvette, France}

\email{kreimer@@ihes.fr}

%\subjclass{17B65, 17B70, 16W30, 81T18, 81T15}

\date{}
\maketitle

\begin{abstract} We prove that the insertion-elimination Lie algebra of Feynman graphs in the ladder case has a natural
interpretation in terms of a certain algebra of infinite
dimensional matrices. We study some aspects of its representation
theory and we will discuss some relations with the representation of the
Heisenberg algebra.
\end{abstract}

\vspace{1cm}

\section{Introduction}

In the last few years perturbative QFT had been shown to have a
rich algebraic structure \cite{Kr} and deep relations with
apparently unrelated sectors of mathematics like non-commutative
geometry and Riemann-Hilbert like problems\cite{C-K 1, C-K 2}.
Such extraordinary relations can be resumed, to some extent, by
the existence of a commutative, non co-commutative Hopf algebra
$\cal H$ defined on the set of Feynman diagrams. In what follows
we will discuss some first relations of perturbative QFT with
the representation theory of Lie algebras. In fact the Hopf
algebra $\cal H$ is isomorphic via the celebrated Milnor-Moore
theorem \cite{M-M} to the dual of the universal enveloping algebra
of a Lie algebra $\cal L$. The Lie algebra $\cal L$ can be
realized as derivations of the Hopf algebra $\cal H$ and it has
been shown in \cite{C-K} to have two natural representations whose
action on $\cal H$ is given by elimination and insertion
operators. Here, the goal of the authors is to study the simplest
case, i.e the case of the cocommutative sub Hopf algebra of ladder
graphs. In this version the Hopf algebra is reduced to a polynomial algebra freely
generated by infinitely many generators and insertion and
elimination are performed by increasing or decreasing the degree
of the generators of such algebra. Even in this simplified context the Lie algebra
introduced in \cite{C-K} has non-trivial features.

A further motivation for the study of this Lie algebra comes from
\cite{BK}. There, a sub Hopf algebra of graphs was studied which
is non-cocommutative and contains the one considered here. The
Dyson-Schwinger equation for that Hopf algebra was solved by
turning the renormalization group into a propagator-coupling
duality. The crucial ingredient was the understanding how the
elimination  Lie algebra (called "befooting" in \cite{BK})
commutes with the generator of the equation of motion, the
Hochschild one-cocycle $B_+$, and how it commutes with the
derivation $S\star Y=m\circ(S\otimes Y)\circ\Delta$ (terminology
as in that paper).

We can easily rewrite the solution of any DSE in the form
$$X=r+\sum_{res(\Gamma)=r}g^{\vert\Gamma\vert}Z_{\Gamma,r}(r),$$
where the residue $r$ specifies the external quantum numbers of
the graph and the sum is over all graphs with external legs
specified by these quantum numbers. Following \cite{BK}, we then
need to understand the commutators
$[Z_{[r,\Gamma]},Z_{[\Gamma,r]}]$ and $[Z_{[r,\Gamma]},S\star Y]$
to solve the DSE in accordance with the RG.

We will express below the derivation $S\star Y$ acting on the Hopf
algebra (of ladder graphs) in terms of generators of the insertion and elimination
Lie algebra, and thus show that the methods of \cite{BK} can be
formulated entirely in that Lie algebra.

\section{The Lie algebra ${\cal L}_{L}$}

Let us start recalling the following theorem where we refere for notation
and symbols to \cite{C-K}:

\begin{thm}\label{t0}
For all 1-PI graphs $\Gamma_{i}$, s.t. {\bf res}$(\Gamma_1)$={\bf res}$(\Gamma_2)$={\bf res}$(\Gamma_3)$={\bf res}$(\Gamma_4)$, the bracket

\begin{eqnarray}
\big[Z_{[\Gamma_{1},\Gamma_{2}]},Z_{[\Gamma_{3},\Gamma_{4}]}\big]
& = &
Z_{\overline{[Z_{[\Gamma_{1},\Gamma_{2}]}\times\delta_{\Gamma_3}},\Gamma_{4}]}-
Z_{[\Gamma_{3},\overline{Z_{[\Gamma_{2},\Gamma_{1}]}\times\delta_{\Gamma_4}}]}
-Z_{\overline{[Z_{[\Gamma_{3},\Gamma_{4}]}\times\delta_{\Gamma_1}},\Gamma_{2}]}\nonumber\\
 & &
+
Z_{[\Gamma_{1},\overline{Z_{[\Gamma_{4},\Gamma_{3}]}\times\delta_{\Gamma_2}}]}-
\delta_{\Gamma_2,\Gamma_3}Z_{[\Gamma_1,\Gamma_4]}+\delta_{\Gamma_1,\Gamma_4}Z_{[\Gamma_3,\Gamma_2]},
\label{00}
\end{eqnarray}
defines a Lie algebra of derivations acting on the Hopf algebra ${\cal H}_{FG}$ via:
$$Z_{[\Gamma_1,\Gamma_2]}\times\delta_{X}=\sum_{I}\langle Z^{+}_{\Gamma_2},\delta_{X'_{(i)}} \rangle\delta_{X''_{(i)}\ast_{G_i}\Gamma_1},$$
where the $G_i$ are normalized gluing data.
\end{thm}

We need to translate the formula (\ref{00}) to the ladder case.
We
start with the following:

\begin{defi}
${\cal L}_{L}= \hspace{5pt}span_{\bC}\hspace{5pt}\{Z_{n,m}; n,m\geq 0\}$

\end{defi}
This is obvious if we identify the $n$-loop ladder graph
$\Gamma_n$ in the Hopf algebra with a Hopf algebra element
$\delta_n$, $n$ being a non-negative integer, and
$Z_{[\Gamma_n,\Gamma_m]}$ with $Z_{n,m}$. Here, any subclass of
$n$-loop graphs $\Gamma_n$ such that
$$\Delta(\Gamma_n)=\Gamma_n\otimes 1+1\otimes
\Gamma_n+\sum_{j=1}^{n-1}\Gamma_j\otimes\Gamma_{n-j}$$ can serve
as an example of a Hopf sub algebra of ladder graphs.

Then, the theorem \ref{t0} becomes:
\begin{thm} ${\cal L}_{L}$ is a Lie algebra with commutator given by the following
formula: \begin{eqnarray} \big[Z_{n,m}, Z_{l,s}\big] & = &
\Theta(l-m)Z_{l-m+n, s}-\Theta(s-n)Z_{l,s-n+m}\nonumber\\ & & -
\Theta(n-s)Z_{n-s+l,m}+\Theta(m-l)Z_{n,m-l+s}\nonumber\\ & &
-\delta_{m,l}Z_{n,s}+\delta_{n,s}Z_{l,m},\label{e5}
\end{eqnarray}
where:\\
\begin{equation}
\begin{cases}
\Theta(l-m) &=0 \hspace{5pt} \text{if $l<m$}, \\
\Theta(l-m) &=1 \hspace{5pt} \text{if $l\geq m$}
\end{cases}
\end{equation}
and
\begin{equation}
\begin{cases}
\delta_{n,m} &=1 \hspace{5pt} \text{if $m=n$},\\
\delta_{n,m} &=0 \hspace{5pt} \text{if $n\neq m$}.
\end{cases}
\end{equation}
\end{thm}

Let us now introduce some natural module for the Lie algebra ${\cal L}_{L}$.

\begin{defi}\label {d1}
$${\cal S}=\bigoplus_{n\geq 0}{\bC}t_n={\bC}[t_0,t_1,t_2,t_3.....].$$
We will assign a degree equal to $k$ to the generator $t_{k}$ for each $k\geq 0$.\\
\end{defi}

${\cal L}_{L}$ acts on $\cal S$ via the following:

\begin{gather}
Z_{n,m}t_k=0\hspace{5pt}if\hspace{5pt} m>k,\notag\\
Z_{n,m}t_k=t_{k-m+n}\hspace{5pt}if\hspace{5pt} m\leq k\label{e0}.
\end{gather}

\begin{rk}\label {r1}
$\cal S$ is a polynomial algebra with infinite many generators graded as in definition \ref{d1}. The product of two polynomials $x_i$, $x_j$
will be denoted by concatenation:
$$x_i\otimes x_j\longrightarrow x_ix_j;$$
in particular $deg(t_it_j)=i+j$ and $y$ is a unit for $\cal S$ with respect to this product if and only if it is a scalar.
 Beside this algebraic structure we need to define on $\cal S$ the following product:
$$\star:{\cal S}\otimes {\cal S}\longrightarrow {\cal S}$$
$$\star(t_{n}\otimes t_m)\mapsto t_{n+m}.$$
With respect to this product $\cal S$ becomes a standard polynomial algebra in one generator with the usual grading:
$$deg(t_k)=k,\hspace{20pt}deg(t_k\star t_l)=deg(t_{k+l})=l+k.$$
In particular $y\in\cal S$ is a unit with respect to the $\star$-product if and only if $y\in \bC t_0$.
\end{rk}

In what follows we will call $\cal S$ the standard representation
for the Lie algebra ${\cal L}_{L}$.

We recall that a Lie algebra $\frak g$ is called $\bZ$-graded if:
$${\frak g}=\bigoplus_{n\in \bZ}{\frak g}_n, \hspace{30pt} [{\frak g}_i,{\frak g}_j]\subseteq {\frak g}_{i+j}.$$
In the decomposition ${\frak g}=\bigoplus_{n\in \bZ}{\frak g}_n$, the components ${\frak g}_{n}$ are called homogeneous of degree equal to $n$.

\begin{prop}
The Lie algebra ${\cal L}_{L}$ is $\bZ$-graded,
$${\cal L}_{L}=\bigoplus_{n\in {\bZ}}l_n,$$
where the homogeneous components of degree $n$ are:
$$l_n=\hspace{5pt}span_{\bC}\{Z_{k,m};\hspace{5pt} k-m=n\}.$$
\end{prop}
\prf We need to prove that if $Z_{n,m}\in l_i$ and $Z_{l,s}\in l_j$ than:
$$[Z_{n,m},Z_{l,s}]\in l_{i+j}.$$
This follows by direct computation using (\ref{e5}).
\endprf

From the $\bZ$-grading it follows that ${\cal L}_{L}$ has the
following decomposition:

$${\cal L}_{L}= L^{+}\oplus L^{0}\oplus L^{-};$$

where $L^{+}=\oplus_{n>0}l_n$, $L^{-}=\oplus_{n<0}l_n$ and $L^{0}=l_0$.

We have that:
\begin{cor}
$L^{+}$, $L^{-}$ and $L^{0}$ are sub Lie algebras of ${\cal L}_{L}$. Moreover $L^{0}$ is commutative and fulfills the following:
$$[L^{+},L^{0}]\subseteq L^{+},\hspace{20pt}[L^{-},L^{0}]\subseteq L^{-}.$$
\end{cor}
\prf It follows by direct computation using (\ref{e5}).
\endprf

We thus have the following:

\begin{prop}\label{l2} In terms of the standard representation we have:\\
given $n,m\geq 0$, $Z_{n,m}$ is a matrix in which the only non-zero
entries are all equal to one and are located  on the $n-m$ lower diagonal if $n>m$, or on the $m-n$ upper
diagonal if $m>n$. More precisely given $n>m\geq 0$ and $k>0$:\\
a) $Z_{n-m,0}$ is a matrix in which the entries of the $n-m$
(lower) diagonal are all equal to 1;\\
b) $Z_{n-m+k,k}$ is a matrix having zeros on the first $k-1$ entries
of the $n-m$ (lower) diagonal and all the other entries of such
diagonal equal to one.\\
a) and b) hold for the matrices $Z_{0,n-m}$ and $Z_{k, n-m+k}$  substituting, in the previous statements, lower with upper.
\end{prop}
\prf
It follows directly from (\ref{e0}).
\endprf
The following claim is now evident:

\begin{lemma}
Given $n,m$ as in the lemma (\ref{l2}), we have that:
if $n-m>0$ $Z_{n,m}\in L^{+}$, if $n-m<0$ $Z_{n,m}\in L^{-}$ and if n=m $Z_{n,m}\in L^{0}$.\\
\end{lemma}

We recall that given a graded Lie algebra ${\frak g}=\bigoplus_{n\in\bZ}{\frak g}_n$, a vector space
$V$ is a graded $\frak g$-module if $V=\bigoplus_{n\in\bZ}V_n$ and $g_jV_i\subseteq V_{i+j}$, $\forall i,j\in\bZ$.
The $\frak g$-module $V$ will be said to be of finite type if ${\text {dim}}\hspace{4pt}V_j<\infty$, $\forall j\in\bZ$.
So we have:
\begin{prop}
$\cal S$ is a graded ${\cal L}_{L}$-module of finite type.
\end{prop}
\prf It follows from direct computation; given any element $Z_{n,m}\in l_i$ with $i=n-m$:
$$Z_{n,m}t_k=t_{k-m+n}\in l_{k+i}$$
from (\ref{e0}).
\endprf

We also recall that a highest weight (h.w) $\frak g$-module of highest weight $\alpha\in{\frak g}_{0}^{\ast}$ is a $\bZ$-graded $\frak g$-module
$V(\alpha)=\bigoplus_{n\in{\bZ}_{\geq 0}}V_n$ such that the following properties hold:\\
a) $V_{0}={\bC}w_{\alpha}$, where $w_{\alpha}$ is a vector not equal to zero;\\
b) $hw_{\alpha}=\alpha (h)w_{\alpha}$, for all $h\in {\frak g}_{0}$;\\
c) ${\frak g}_{-}w_{\alpha}=0$;\\
d) ${\cal U}({\frak g}_{+})w_{\alpha}=V(\alpha)$,\\
where ${\frak g}_{+}=\bigoplus_{n>0}{\frak g}_{n}$, ${\frak g}_{-}=\bigoplus_{n<0}{\frak g}_{n}$ and with ${\cal U}{(\frak g)}$ we indicated the universal enveloping algebra of $\frak g$.
The vector $w_{\alpha}$ is called highest weight vector (h.w.v) and any vector $v$ such that ${\frak g}_{-}v=0$ is called singular vector. Moreover
we have that the module $V(\alpha)$ is irreducible if and only if every singular vector is a multiple of the h.w.vector.
Now we can state the following:
\begin{prop}
$\cal S$ is an irreducible h.w. module for the Lie algebra ${\cal L}_{L}$ with h.w. vector $w_{\alpha}=t_{0}$.
\end{prop}
\prf It follows directly by the definitions.

\endprf

\section{Classical Lie algebras and the Lie algebras ${\cal L}_{L}$.}

In what follows we will investigate some relations of the
Lie algebra ${\cal L}_{L}$ with some (infinite dimensional) classical Lie algebras. In particular
we will give a complete description of ${\frak sl}_{+}(\infty)$ in terms of ${\cal L}_{L}$.
Let's start with the following:

\begin{defi}\cite{K, K-R}

$${\frak gl}(\infty)=\{ E_{i,j}:\hspace{5pt}i,j\in {\bZ}: [E_{i,j}, E_{n,m}]=E_{i,m}\delta_{j,n}-E_{n,j}\delta_{m,i}\};$$

and

$${\frak gl}_{+}(\infty)=\{ E_{i,j}:\hspace{5pt}i,j\geq 0: [E_{i,j}, E_{n,m}]=E_{i,m}\delta_{j,n}-E_{n,j}\delta_{m,i}\}.$$

\end{defi}

Now we can state the following:

\begin{prop}
We have an embedding of Lie algebras:
\begin{equation}
\phi:{\frak gl}_{+}(\infty)\longrightarrow{\cal L}_{L}. \label{i}
\end{equation}
\end{prop}
\prf  Define $\phi(E_{i,j})=Z_{i,j}-Z_{i+1,j+1}$ for each $i,j\geq 0$.\\
Now it suffices to show that $\phi$ is morphism of Lie algebras, i.e that:
$$[Z_{i,j}-Z_{i+1,j+1}, Z_{n,m}-Z_{n+1,m+1}]=(Z_{i,m}-Z_{i+1,m+1})\delta_{j,n}-(Z_{n,j}-Z_{n+1,j+1})\delta_{m,i}.$$
This follows directly from the definitions.
\endprf

\begin{defi}\cite{K, K-R}
$${\frak sl}_{+}(\infty)=\{A\in{\frak gl}_{+}(\infty):\hspace{5pt}trA=0\}$$
\end{defi}
We recall now that a set of generators for ${\frak sl}_{+}(\infty)$ is given by:
\begin{equation}
\begin{cases}
E_{i,j}&  i<j,\\
E_{i,j}&  i>j, \\
E_{i,i}-E_{l,l}& \text{with $i\neq l$}.\label{e1}
\end{cases}
\end{equation}

From (\ref{e1}) we get the following (standard) triangular
decomposition:
\begin{equation}
{\frak sl}_{+}(\infty)={\frak n}_{+}\oplus{\frak h}\oplus{\frak n}_{-};\label{e3}
\end{equation}
where

${\frak n}_{+}=span_{\bC}\{E_{i,j};\hspace{5pt}j>i\}$; ${\frak n}_{-}=span_{\bC}\{E_{i,j};\hspace{5pt}i>j\}$;
${\frak h}=span_{\bC}\{E_{i,i}-E_{l,l};\hspace{5pt}i\neq l\}$.\\
Let us now introduce Chevalley's generators and co-roots for the Lie algebra ${\frak sl}_{+}(\infty)$ (\cite {K}):
The Chevalley's generators are:

$$e_i=E_{i,i+1}\hspace{20pt} f_{i}=E_{i+1,i}\hspace{8pt}\forall i\in{\bZ}_{\geq 0}$$

and

$$\check\Pi=\{E_{i,i}-E_{i+1,i+1}; i\in{\bZ}_{\geq 0}\}$$

is the set of simple co-roots.

We can now write the Chevalley's \cite{K} generators and a co-root system for ${\frak
sl}_{+}(\infty)$ in terms of the generators $Z_{n,m}$ of the Lie
algebra ${\cal L}_{L}$:
\begin{lemma}
The Chevalley generators for ${\frak sl}_{+}(\infty)$, $\{e_i,
f_i;\hspace{5pt}i\geq 0\}$ can be written in
terms of the generators $Z_{n,m}$ in the following way:
\begin{equation}
\begin{cases}
f_i &= Z_{i+1,i}-Z_{i+2,i+1},\\
e_i &= Z_{i,i+1}-Z_{i+1,i+2}.
\end{cases}
\end{equation}
If we define:
$$\check{\alpha}_i = Z_{i,i}-2Z_{i+1,i+1}+Z_{i+2,i+2};$$
then $\{\check{\alpha_i};\hspace{5pt}i\geq 0\}$ is a system of positive simple co-roots for the Lie algebra ${\frak sl}_{+}(\infty)$.
\end{lemma}
\prf It follows from the definition of Chevalley generators, positive simple roots and
from the embedding of ${\frak gl}_{+}(\infty)$ in ${\cal L}_{L}$ defined in (\ref{i}).
\endprf

\begin{rk}\cite{K}
The root system for ${\frak sl}_{+}(\infty)$ is described in the
following way; we remember that we have a triangular decomposition
for ${\frak gl}_{+}(\infty)$ analogous to (\ref{e3}):
$${\frak gl}_{+}(\infty)={\frak n}_{+}\oplus{\frak h}\oplus{\frak n}_{-},$$
where now:
$${\frak h}=span_{\bC}\{E_{i,i};\hspace{5pt} i\geq 0\}.$$
Define now:
$${\frak h}^{\ast}=\{{\bC}-\text{linear functions on}\hspace{5pt}{\frak h}\}=span_{\bC}\{\epsilon_i;\hspace{5pt} i\geq 0\},$$
where
$$\epsilon_i(E_{j,j})=\delta_{i,j},$$
or in terms of $Z_{n,m}$ generators:
$$\epsilon_i(Z_{j,j}-Z_{j+1,j+1})=\delta_{i,j}.$$
From this description follows immediately:\\ The root system for
${\frak sl}_{+}(\infty)$ is given by:
\begin{itemize}
\item[i)] $\Delta=span_{\bC}\{\epsilon_i-\epsilon_j;\hspace{5pt}i\neq j\}$;

\item[ii)] $\Delta^{+}=span_{\bC}\{\epsilon_i-\epsilon_j;\hspace{5pt}i,j\geq 0;\hspace{5pt}i<j \}$;

\item[iii)] $\Pi=span_{\bC}\{\epsilon_i-\epsilon_{i+1};\hspace{5pt}i\neq 0\}$;

\end{itemize}
where $\Delta$, $\Delta^{+}$ and $\Pi$ are
respectively called the set of roots, the set of positive roots and the set of simple positive roots.

\end{rk}

\subsection{The Chevalley's involution on ${\cal L}_{L}$.}

We define now an involution on ${\cal L}_{L}$ whose restriction to ${\frak sl}_{+}(\infty)$ gives us the Chevalley's involution.
Let us start with:
\begin{defi}
$$C:{\cal L}\longrightarrow{\cal L}$$
$$Z_{n,m}\longmapsto -Z_{m,n}$$
for each $n,m\geq 0$.
\end{defi}
We have the following:
\begin{prop}
$C$ is an homomorphism of Lie algebras, i.e:
\begin{equation}
C([Z_{n,m},Z_{l,s}])=[C(Z_{n,m}),C(Z_{l,s})],\label{e4}
\end{equation}
$\forall\hspace{5pt} n,m, l, s\geq 0$.
\end{prop}
\prf
It follows applying the (\ref{e5}) to both sides of (\ref{e4}):\\

RHS=$\big[Z_{m,n},Z_{s,l}\big]=\Theta(s-n)Z_{s-n+m, l}-\Theta(l-m)Z_{s,l-m+n}-$\\
\begin{equation}
\Theta(m-l)Z_{m-l+s,n}+\Theta(n-s)Z_{m,n-s+l}
-\delta_{n,s}Z_{m,l}+\delta_{m,l}Z_{s,n},
\end{equation}

LHS=$C\big(\big[Z_{n,m},Z_{l,s}\big]\big)=-\Theta(l-m)Z_{s,l-m+n}+\Theta(s-n)Z_{s-n+m,l}+$\\
\begin{equation}
\Theta(n-s)Z_{m,n-s+l}-\Theta(m-l)Z_{m-l+s,n}
-\delta_{n,s}Z_{m,l}+\delta_{m,l}Z_{s,n}.
\end{equation}
\endprf

We recall that:
\begin{defi}\cite{K}
Given a Lie algebra $\frak g$, a set of its Chevalley's
generators $\{f_i,e_i;i\geq 0\}$ and a system of simple positive co-roots $\{\check{\alpha}_i,\hspace{5pt}i\geq 0\}$, the
map:
$$\omega:{\frak g}\longrightarrow {\frak g}$$
defined by:
$$\omega(f_i)=-e_i,\hspace{8pt}\omega(e_i)=-f_i,\hspace{8pt}\omega(\check{\alpha}_i)=-\check{\alpha}_i,$$
is called Chevalley's involution.
\end{defi}

\begin{prop}
The restriction of the map $C$ to ${\frak sl}_{+}(\infty)$ is the
Chevalley's involution.
\end{prop}
\prf It is obvious from the definitions:
$$C(f_i)=C(Z_{i+1,i}-Z_{i+2,i+1})=-Z_{i,i+1}+Z_{i+1,i+2}=e_{i},$$
$$C(\check{\alpha}_{i})=C(Z_{i,i}-2Z_{i+1,i+1}+Z_{i+2,i+2})=-Z_{i,i}+2Z_{i+1,i+1}-Z_{i+2,i+2}={\check{\alpha}}_{i}.$$
\endprf

\section {$l^{+},\hspace{5pt}l^{-}$ and their central extensions.}

In this section we are going to consider the relation between the Lie algebra ${\cal L}_L$ and the Heisenberg algebra.
We will introduce two sub-algebras of ${\cal L}_{L}$ that play the role of "shift" operators for the
standard module $\cal S$.
Let us start with the following:
\begin{defi}
Define:
\begin{equation}
l^{+}=span_{\bC}\{Z_{n,0};\hspace{5pt} n\geq 0\}
\end{equation}
and
\begin{equation}
l^{-}=span_{\bC}\{Z_{0,n};\hspace{5pt} n\geq 0\}.
\end{equation}
\end{defi}

We have:
\begin{lemma}\label{l1}
1)$$\big[l^{+},l^{+}\big]=0,\hspace{15pt}\big[l^{-},l^{-}\big]=0;$$
2)\hspace{6pt} $l^{+}$ acts as algebra of shift operators on $\cal S$ and $l^{-}$ as algebra of quasi-shift operators, i.e
$\forall n, k>0$ $Z_{0,n}$ will eliminate any $t_k\in\cal S$ if $k<n$ and will shift such element to the element $t_{k-n}$ if $n<k$.
\end{lemma}
\prf
a) follows applying (\ref{e5}) to elements in $l^{+}$ and $l^{-}$.\\
b) follows trivially from the definition of the ${\cal L}_{L}$ action on $\cal S$.
\endprf

Note that the commutativity of these two algebras stems from the
cocommutativity of the ladder Hopf algebra ${\cal H}_L$. The
general case provides non-commutative Lie algebras for insertion
as well as elimination \cite{C-K}.

 Let us
now define two other commutative Lie algebras.
\begin{defi}
Let $\{Z_{n}^{+}; n\in {\bZ}\}$ and $\{Z_{n}^{-}; n\in {\bZ}\}$
two sets of symbols. Let us define:
$$\check{l^{+}}=span_{\bC}\{Z_{n}^{+};\hspace{5pt}n\in {\bZ}\}$$
and
$$\check{l^{-}}=span_{\bC}\{Z_{n}^{-};\hspace{5pt}n\in {\bZ}\}$$
and let us introduce the canonical isomorphism:
$$d:\check{l^{+}}\longrightarrow\check{l^{-}}$$
defined on the generators by:
$$d(Z_n^{+})=Z_n^{-}\hspace{15pt} n\in {\bZ}$$
\end{defi}
Now define the following maps:

$${\frak a}^{+}:l^{+}\longrightarrow\check{l^{+}}$$
$$Z_{2n,0}\mapsto Z_{n}^{+}\hspace{5pt}\forall\hspace{5pt}n\geq 0$$
$$Z_{2n-1,0}\mapsto Z_{-n}^{+}\hspace{5pt}\forall\hspace{5pt}n\geq 1$$
and

$${\frak a}^{-}:l^{-}\longrightarrow\check{l^{-}}$$
$$Z_{0,2n}\mapsto Z_{n}^{-}\hspace{5pt}\forall\hspace{5pt}n\geq 0$$
$$Z_{0,2n-1}\mapsto Z_{-n}^{-}\hspace{5pt}\forall\hspace{5pt}n\geq 1$$

Now we have the following:
\begin{lemma}
${\frak a}^{+}$ and ${\frak a}^{-}$ are isomorphisms of Lie algebras compatible with the involution $C$.
\end{lemma}
\prf
It follows trivially from the definitions of $C$, $d$ and from the commutativity of the Lie algebras $l^{\pm}$ and $\check{l^{\pm}}$.
\endprf

It is a well known result that:
\begin{prop}\cite{K, K-R}
$$H^{2}(\check{l^{\pm}},{\bC})\neq 0$$
\end{prop}
\prf
Let's focus on $\check{l^{+}}$ case (the case $\check{l^{-}}$ is completely analogous). Define the following bilinear map:
$$c^{+}:\check{l^{+}}\otimes\check{l^{+}}\longrightarrow{\bC}$$
$$(Z_{n}^{+}\otimes Z_{m}^{+})\mapsto n\delta_{n,-m}.$$
$c^{+}$ clearly satisfies the following:
$$c^{+}(Z_{n}^{+},Z_{m}^{+})=-c^{+}(Z_{m}^{+},Z_{n}^{+})\hspace{5pt}{\text and}$$
$$c^{+}\big(\big[Z_{n}^{+},Z_{m}^{+}\big],Z_{k}^{+}\big)+\hspace{5pt}\text{cyclic permutations}=0;$$
i.e $c^{+}$ is a non trivial two co-cycle.
\endprf

\begin{rk}\cite{K, K-R}
The co-cycle gives rise to a central extension:
\begin{equation}
0@>>> {\bC}{\cal C} @>>> {\cal H}^{+} @>p>> \check{l^{+}} @>>> 0,
\end{equation}
where ${\cal H}^{+}=span_{\bC}\{Z_{n}^{+}, {\cal C}\}_{n\in
{\bZ}}$ and relations:
$$\big[Z^+_{n},Z^+_{m}\big]=n\delta_{n,-m}{\cal C}\hspace{15pt} \big[Z^+_{n},{\cal C}\big]=0.$$

${\cal H}^{+}$ is called Heisenberg Lie algebra; the "$-$" case is completely analogous: the co-cycle $c^{-}$ will give us the central extension
${\cal H}^{-}$ with ${\cal H}^{-}\cong {\cal H}^{+}$.

Let us now consider a $\cal H$-module V (in what follows we will drop the $\pm$ sign) on which $Z_0$ acts as multiplication
operator. Let us define the following operators:

$$L_{0}=\frac{{\mu}^2+{\lambda}^2}{2}+\sum_{n>0}Z_{-n}Z_n,$$
and
$$L_n=\frac{1}{2}\sum_{j\in\bZ}Z_{-j}Z_{j+n}+i\lambda nZ_n\hspace{5pt}\forall n\in{\bZ}^{\ast}\hspace{5pt}\lambda\in{\bC}.$$
Then we have that the $L_n$'s defined above fulfill the Virasoro commutation relations:
$$\big[L_n,L_m\big]=(n-m)L_{n+m}+\delta_{n,-m}\frac{(n^3-n)}{12}(1+12\lambda^2)Id_{V},$$
where $\lambda\in\bC$ gives the central charge of the Virasoro algebra.

\end{rk}

It is now time to turn to the derivation $S\star Y$, expressed in generators $Z_{0,n},Z_{m,0}$ and their commutators.

\section{The derivation $S\star Y$}

Let us now turn to an expression of the derivation $S\star Y$ in
terms of the generators of the insertion elimination Lie algebra.
To that end, let $\Gamma_m$ correspond to the $m$-loop ladder
graph as a generator in the sub Hopf algebra ${\cal H}_L$ of
ladder graphs.

Define a derivation $D_1$ on ${\cal H}_L$ by
$$D_1(\Gamma_m)=\sum_{n=0}^\infty \Gamma_n Z_{0,n}(\Gamma_m) $$
which takes care of $m\circ \Delta$: for any two characters
$\phi_1,\phi_2$ on ${\cal H}_L$ we have
$$\phi_1\star\phi_2(\Gamma_m)=m\circ (\phi_1\otimes \phi_2)\circ\Delta(\Gamma_m)=\sum_{n=0}^\infty \phi_1(\Gamma_n)\phi_2(Z_{0,n}(\Gamma_m)).$$
The degree operator $Y$ can be described by a derivation $D_2$,
$$D_2(\Gamma_m)=\sum_{k=1}^\infty Z_{k,k}(\Gamma_m)=m\Gamma_m $$
and $S(\Gamma_m)$ iteratively by a derivation $D_3$,
$$D_3(\Gamma_m)=-Z_{0,0}(\Gamma_m)-\sum_{n=0}^\infty
D_3(\Gamma_n)Z_{1,n+1}(\Gamma_m).
$$
Note that $Z_{k,k}=[Z_{k,0}, Z_{0,k}]+Z_{0,0}$, similarly $Z_{1,n+1}$ involves the commutator $[Z_{1,0},Z_{0,n+1}]$.\\
Now we get $S\star Y$ from the definitions.

\begin{prop}
$S\star Y$ is the derivation that is uniquely given on the linear generators $\Gamma_i$ as:
\begin{equation}
S\star Y(\Gamma_m)=\sum_{n=0}^\infty D_3(\Gamma_n)D_2(Z_{0,n}(\Gamma_m)).\label{S1}
\end{equation}
\end{prop}

\begin{rk}
Note that (\ref{S1}) holds as the coproduct is linear on generators on the lhs in the ladder case,
and on the linear subspace of generators the antipode can indeed be
decribed by a derivation.
\end{rk}

\section{The standard module $\Lambda$.}

In this section we are going to address the following question: what is the equivalent of the standard
module $\cal S$ for the Lie algebras $\check{l^{\pm}}$?
Recall (Remark \ref{r1}) that on ${\cal S}=\bigoplus_{n\geq 0}{\bC}t_{n}$ one can naturally define the product:
$$\star:{\cal S}\otimes {\cal S}\longrightarrow {\cal S}$$
$$\star(t_{n}\otimes t_m)\mapsto t_{n+m}$$
and that with respect to this product $\cal S$ is isomorphic to a polynomial algebra having only one generator.
We also remember that the action of $l^{\pm}$ is given in \ref{l1} (2). To define the
standard module for $\check{l^{\pm}}$, we introduce the following
notation:
$$o(n)\doteq -exp(n-\frac{1}{2}),\hspace{5pt}n>0;$$
$$e(n)\doteq exp(n),\hspace{5pt}n\geq 0.$$
\begin{defi}
Let us start defining the vector space:
$$\Lambda=span_{\bC}\{\alpha(o(n)), \alpha(e(n))\}$$
$\Lambda$ is a unital algebra with the following product:
$$\bullet :\Lambda\otimes\Lambda\longrightarrow\Lambda$$
$$\alpha(\xi(n))\bullet\hspace{3pt} \alpha(\xi(m))=\alpha(\xi(n)\xi(m));$$
where $\xi(k)$can be either $o(k)$ or $e(k)$. The unit is given by $\alpha=\alpha(e(0))$.
\end{defi}

Let us now define:

$$\phi:{\cal S}\longrightarrow\Lambda$$
by the following:
$$\phi(t_{2k})=\alpha(e(k)),\hspace{5pt}k\geq 0\hspace{5pt}\phi (t_{2k-1})=\alpha(o(k)),\hspace{5pt}k>0.$$
\begin{prop}
The map $\phi$ is an isomorphism of $\bC$-algebras.
\end{prop}
\prf
We need only to check that $\phi$ is a morphism of algebras, for example:
$$\phi(t_{2n-1}\star t_{2m-1})=\phi(t_{2(n+m-1)})=$$
$$\alpha(e(n+m-1))=\alpha(o(n)o(m))=$$
$$\alpha(o(n))\bullet\hspace{3pt}\alpha(o(m))=\phi(t_{2n-1})\bullet\hspace{3pt}\phi(t_{2m-1}).$$
\endprf

We have now the following:

\begin{prop}
$\Lambda$ is a $\check{l^{\pm}}$ module.
\end{prop}
\prf
It suffices to define:
$$\lambda^{\pm}:\check{l^{\pm}}\longrightarrow End(\Lambda),$$
as multiplication operators since $\check{l^{\pm}}$ are commutative Lie algebras and $\Lambda$ is a commutative algebra.
Let us define:
$$\lambda^{+}(Z_{n}^{+})(\alpha(\xi(k)))\doteq\alpha(e(n))\bullet\hspace{3pt}\alpha(\xi(k)),$$
$$\lambda^{+}(Z_{-n}^{+})(\alpha(\xi(k)))\doteq\alpha(o(n))\bullet\hspace{3pt}\alpha(\xi(k))$$
and:
$$\lambda^{-}(Z_{n}^{-})(\alpha(\xi(k)))\doteq\alpha({\tilde e}(n))\bullet\hspace{3pt}\alpha(\xi(k)),\hspace{3pt}if\hspace{3pt}k-n\geq 0$$
and $\lambda^{-}(Z_{n}^{-})(\alpha(\xi(k)))=0$ otherwise;
$$\lambda^{-}(Z_{-n}^{-})(\alpha(\xi(k)))\doteq\alpha({\tilde o}(n))\bullet\hspace{3pt}\alpha(\xi(k)),\hspace{3pt}if\hspace{3pt}k-n\geq 0,$$
and $\lambda^{-}(Z_{-n}^{-})(\alpha(\xi(k)))=0$ otherwise,
where ${\tilde e}(n)=e(-n)=exp(-n)$\\ while ${\tilde o}(n)=-exp(-n+\frac{1}{2})$.
\endprf

We can now state and prove the main result of this section:

\begin{thm}
The following diagrams commute:
\[
\begin{array}{ccc}

S                                    &

\stackrel{Z_{2n,0}}{\longrightarrow}                   &

S                                            \\
\Big\downarrow\vcenter{%
\rlap{$\phi$}}                               & &
\Big\downarrow\vcenter{%
\rlap{$\phi$}}                               \\
\Lambda                              &

\stackrel{Z^{+}_{n}}{\longrightarrow}                   &

\Lambda

\end{array}
\]

\[
\begin{array}{ccc}

S                                    &

\stackrel{Z_{2n-1,0}}{\longrightarrow}                   &

S                                            \\
\Big\downarrow\vcenter{%
\rlap{$\phi$}}                               & &
\Big\downarrow\vcenter{%
\rlap{$\phi$}}                               \\
\Lambda                              &

\stackrel{Z^{+}_{-n}}{\longrightarrow}                   &

\Lambda

\end{array}
\]
The same statement holds for:
\[
\begin{array}{ccc}

S                                    &

\stackrel{Z_{0,2n}}{\longrightarrow}                   &

S                                            \\
\Big\downarrow\vcenter{%
\rlap{$\phi$}}                               & &
\Big\downarrow\vcenter{%
\rlap{$\phi$}}                               \\
\Lambda                              &

\stackrel{Z^{-}_{n}}{\longrightarrow}                   &

\Lambda

\end{array}
\]
\[
\begin{array}{ccc}

S                                    &

\stackrel{Z_{0,2n-1}}{\longrightarrow}                   &

S                                            \\
\Big\downarrow\vcenter{%
\rlap{$\phi$}}                               & &
\Big\downarrow\vcenter{%
\rlap{$\phi$}}                               \\
\Lambda                              &

\stackrel{Z^{-}_{-n}}{\longrightarrow}                   &

\Lambda

\end{array}
\]
\end{thm}
\prf
It follows by inspection. Let us consider for example the fourth diagram:
$$\phi(Z_{0,2n-1}(t_{2k-1}))=\phi(t_{2(k-n)})=\alpha(e(k-n))$$
and
$$Z_{-n}^{-}(\phi(t_{2k-1}))=Z_{-n}^{-}(\alpha(-exp(k-\frac{1}{2})))=\alpha(-exp(-n+\frac{1}{2}))\bullet\hspace{3pt}\alpha(-exp(k-\frac{1}{2}))=$$
$=\alpha(e(k-n))$. Similarly:
$$\phi(Z_{0,2n-1}(t_{2k}))=\phi(t_{2(k-n)+1})=\alpha(-exp(k-n+\frac{1}{2}))=\alpha(o(k-n+1)),$$
and
$$Z_{-n}^{-}(\phi(t_{2k}))=Z_{-n}^{-}(\alpha(exp(k)))=\alpha(-exp(-n+\frac{1}{2}))\bullet\hspace{3pt}\alpha(exp(k))=$$
$=\alpha(-exp(k-n+\frac{1}{2}))=\alpha(o(k-n+1))$.

\endprf

\section{Conclusions and outlooks}
With this paper we started the study of the Lie algebra introduced
in \cite {C-K}. We considered the simplest case for this
Insertion-Elimination Lie algebra, the one coming from the
sub-Hopf algebra of ladder graphs, giving a description of such a
Lie algebra in terms of (infinite) matrices. We then described the
relations of such Lie algebra with other well known infinite
dimensional Lie algebras like the Heisenberg algebra and ${\frak
gl}(\infty)$.

In forthcoming works we will study the structure and
representation theory of ${\cal L}_L$ and we will consider the
case of of the Insertion-Elimination algebra coming from the
ladder Hopf algebra with additional decorations, which makes the
underlying Hopf algebra non-cocommutative, in contrast to the case
studied here.

\vskip20pt
\paragraph{\bf Acknowledgements.}
I.M.~thanks the IHES for hospitality during a stay in the spring
of '03.  The authors thank Kurusch Ebrahimi Fard for his careful
reading of a preliminary version of the paper and for valuable
discussions.


\begin{thebibliography}{BMP}
\small

\bibitem{BK} D.\ J.\ Broadhurst, D.\ Kreimer
{\it Exact solutions of Dyson-Schwinger equations for iterated
one-loop integrals and propagator-coupling duality,} Nucl.\ Phys.\
B {\bf 600}, 403 (2001) [arXiv:hep-th/0012146].

\bibitem{C-K 1} A.\ Connes, D.\ Kreimer {\it Renormalization in quantum field theory and the Riemann Hilbert problem.
I.\ The Hopf algebra structure of graphs and the amin theorem},
Comm. Math. Phys. {\bf 210} (2000), no.1, 249-273.

\bibitem{C-K 2} A.\ Connes, D.\ Kreimer {\it Renormalization in quantum field theory and the Riemann Hilbert problem
. II.\ The $\beta$-function, diffeomorphism and the
reneormalization group}, Comm. Math. Phys. {\bf 216} (2001), no.1,
215-241.

\bibitem{C-K} A.\ Connes, D.\ Kreimer, {\it Insertion and Elimination: the doubly infinite Lie algebra of Feynmann graphs},
Ann. Henri Poincare {\bf 3}, (2002) no. 3, 411-433.

\bibitem{K}  V.\ Kac {\it Infinite Dimensional Lie Algebras}, Cambridge University Press, 3rd Ed. (1994).

\bibitem{K-R} V.\ Kac, A.\ Raina {\it Bombay Lectures on Highest Weight Representation of Infinite Dimensional Lie Algebras},
Advanced Series in Mathematical Physics, Vol 2, World Scientific Pub., 1988.


\bibitem{Kr} D.\ Kreimer {\it On the Hopf algebra structure of perturbative quantum field theory},
Adv. Theor. Math. Phys, {\bf 2}, no. 2, 1998.


\bibitem{M-M} J.\ W.\ Milnor, J.\ C.\ Moore {\it On the Structure of Hopf Algebras},
Ann. of Math, {\bf 81}, no. 2, 1965, 211-264.


\end{thebibliography}
\end{document}